\input amstex.tex

\documentstyle{amsppt}
\magnification \magstephalf
\pagewidth{6.2in}
\pageheight{8.0in}
\NoRunningHeads
\def\makefootline{\baselineskip=49pt \line{\the\footline}}
\define\s{\sigma}
\define\la{\lambda}
\define\al{\alpha}

\define\ep{\varepsilon}

\define\trace{\operatorname{trace}}
\define\supp{\operatorname{supp}}
\define\PDO{\Psi\text{DO}}
\def\today{\ifcase\month\or
    January\or February\or March\or April\or May\or June\or
    July\or August\or September\or October\or November\or December\fi
    \space\number\day, \number\year}

\topmatter
\title
Extremals for Logarithmic Hardy-Littlewood-Sobolev inequalities on
compact manifolds
\endtitle
\author  K. Okikiolu
\endauthor
\thanks  The author was supported  by the National
Science Foundation  \#DMS-0302647.  \endthanks
\abstract Let $M$ be a closed, connected surface and let $\Gamma$ be
a conformal class of metrics on $M$ with each metric normalized to
have area $V$.  For a metric $g\in \Gamma$, denote the area element
by $dV$ and the Laplace-Beltrami operator by $\Delta_g$.  We define
the Robin mass $m(x)$ at the point $x\in M$ to be the value of the
Green's function $G(x,y)$ at $y=x$ after the logarithmic singularity
has been subtracted off.  The regularized trace of $\Delta_g^{-1}$
is then defined by $\trace \Delta^{-1}=\int_M m\,dV$.  (This
essentially agrees with the zeta functional regularization and is
thus a spectral invariant.) Let $\Delta_{S^2,V}$ be the
Laplace-Beltrami operator on the round sphere of volume $V$.  We
show that if there exists $g\in \Gamma$ with $\trace
\Delta_g^{-1}<\trace \Delta_{S^2,V}^{-1}$ then the  minimum of
$\trace \Delta^{-1}$ over $\Gamma$ is attained by a metric in
$\Gamma$ for which the Robin mass is constant.  Otherwise, the
minimum of $\trace \Delta^{-1}$ over $\Gamma$ is equal to $\trace
\Delta_{S^2,V}^{-1}$.  In fact we prove these results in the general
setting where $M$ is an $n$ dimensional closed, connected manifold
and the Laplace-Beltrami operator is replaced by any non-negative
elliptic operator $A$  of degree $n$ which is conformally covariant
in the sense that for the metric $g$ we have $A_{F^{2/n}g}=F^{-1}
A_g$.  In this case the role of $\Delta_{S^2,V}$ is assumed by  the
{\it Paneitz} or {\it GJMS operator} on the round $n$-sphere of
volume $V$. Explicitly these results are logarithmic HLS
inequalities for $(M,g)$. By duality we obtain analogs  of the
Onofri-Beckner theorem.
\endabstract
\endtopmatter

\noindent{\bf Section 1.  Introduction and Results.}
\medskip

Let $S^n$ denote the standard unit sphere in $\Bbb R^{n+1}$  with normalized volume element $d\s$.
In [CL] and [Be] the following endpoint Hardy-Littlewood-Sobolev inequality was established:
\proclaim{\bf Sharp logarithmic Hardy-Littlewood-Sobolev inequality on the sphere}
The inequality
$$
 \frac2{n!}\int_{S^n} F\log F\,d\s\ -\ \int_{S^n}\, F \square^{-1} F\,d\s \ \geq\ 0
\tag 1.1
$$
holds for all functions $F:S^n\to [0,\infty)$ with $ \int_{S^n} F\,d\s=1$,
such that $\int_{S^n} F\log F\,d\s$ is finite.  Moreover equality is attained exactly
when $F$ is the Jacobian of a conformal transformation of $S^n$.
\endproclaim

Here, the operator $\square$ is a natural operator of order $n$ on $S^n$ given by its action on
a spherical harmonic $Y_k$ of degree $k$ on $S^n$ by
$$
\square Y_k\ =\ k(k+1)\dots (k+n-1)Y_k.
$$
When  $n=2$, $\square$ is just the Laplace-Beltrami operator, when $n=4$  it is
the Paneitz operator, and for general even $n$ it
is the GJMS operator.
The inverse
$$
\square^{-1} Y_k\ =\ \cases Y_k/ (k(k+1)\dots (k+n-1)), \qquad & k>0\\ 0 & k=0\endcases
$$
can be written as an integral operator.  Indeed, writing $|x|$ for the Euclidean norm
of $x\in\Bbb R^{n+1}$, we have
$$
\square^{-1} F(x)\ =\ -\frac{2}{(n-1)!}\int_{S^n} \log|x-y|\ (F(y)-\s_F)\, d\s(y),\qquad\qquad\text{ where }
\s_F:=\int_{S^n}F\,d\s.
$$
The inequality (1.1) is dual to the following inequality, see [Be], [CL], {CY], which in
$2$ dimensions is known as Onofri's inequality [On].
\proclaim{\bf Beckner's inequality} If $u$ is in the Sobolev space $L^2_{n/2}(S^n)$  (functions on
$S^n$ with $n/2$ derivatives in $L^2$), then
$$
\frac1{2n!} \int_{S^n} u \square u\,d\s\ + \ \int_{S^n} u\,d\s\ -\ \log\left(\int_{S^n} e^{u}\,d\s\right)\
\geq\ 0.  \tag 1.2
$$
 Moreover equality is attained exactly
when $e^u$ is a multiple of the Jacobian of a conformal transformation of $S^n$.
\endproclaim
For some related inequalities, see [Ad], [Au], [CC], [F], [L], [Mor2], [Mos], [Ok*].

Now the inequalities (1.1) and (1.2) have geometric interpretations.
The best known  is the case of Onofri's inequality, for which the
interpretation involves the zeta-regularized determinant of the
Laplace-Beltrami operator.

 \proclaim{{\bf Polyakov-Ray-Singer
formula}} Let $M$ be a closed surface with metric $g$, area element
$dV$, and let $K$ be the Gaussian curvature of $M$.
$$
\log \det \Delta_{e^u g}\ -\ \log\det \Delta_g\ =\ -\frac1{12\pi}\left(
\frac14\int_{M} |\nabla u|^2 \,dV\ +\ \int_{M}K u\,dV\right)\ +\
\log \int_{M} e^u\,\frac{dV}V. \tag 1.3
$$
\endproclaim

We see that when $g$ is the standard metric  on $S^2$ and the area
of $(M, e^ug)$ equals that of $(M,g)$, then the Polyakov functional
in (1.3) is a multiple of the left hand side of (1.2).  Hence
Onofri's inequality can  be interpreted as saying that among all
metrics conformal to the standard metric on $S^2$ having the same
area, the standard metric attains the maximum value of $\det\Delta$.
See [OPS1]. The analog of the Polyakov-Ray-Singer formula for the
determinant of $\square$ in $4$ dimensions was computed  by Branson
[Br].  On $S^4$, the Beckner functional does occur in this formula,
plus another functional which is also minimized at the standard
metric. For other extremal results concerning determinants, see for
example [Br\O],
 [CQ], [CY], [HZ], [Ok1], [Ok2], [OW],
[OPS2], [R].
\medskip

The inequality (1.1) on the other hand is related to the the regularized trace of $\square^{-1}$.
Indeed, suppose $g$ is the standard metric on $S^n$ and $F$ is a positive smooth function on $S^n$.
Then for the metric $F^{2/n}g$, define the operator $\square$
 by $\square_{F^{2/n}g}=F^{-1}\square_ g$.
The leading order term of
$\square$ agrees with $\Delta^{n/2}$ where $\Delta$ is the Laplace-Beltrami operator.  When $n$ is even,
the differential operator
$\square$ is known as the GJMS operator and is given locally by a geometrically invariant formula.
The following
formula was proved in [Mor1] for the zeta regularization, and then [Mor2], [St1], [St2].
(See also (1.5) and (1.6).)
\proclaim{\bf Conformal change of $\trace\square^{-1}$ on the sphere}
 If $g$ is the standard metric on $S^n$
and  $\int_M F\,d\s=1$, then
$$
\trace \square_{F^{2/n}g}^{-1}\ -\ \trace \square_g^{-1}\ =\ \frac{2}{n!}\int_{S^n} F\log F\,d\s
\ -\ \int_{S^n} F\square^{-1} F\,d\s.
$$
\endproclaim
Morpurgo [Mor1] then interpreted the sharp logarithmic HLS
inequality  as saying that among all metrics conformal to $g$ with
the same volume, the standard metric minimizes $\trace
\square^{-1}$. We also note that Steiner [St1], [St2] gave an
interpretation of the sharp logarithmic HLS inequality as an analog
of the Riemannian positive mass theorem, see [SY], [AH].
\medskip

In this paper, we seek to understand analogs of (1.1), (1.2) on
general closed manifolds. Let $M$ be a smooth compact $n$
dimensional manifold without boundary with a Riemannian metric
$g_0$.  We denote the volume element of $g_0$ by $dV_0$ and the
volume of $(M,g_0)$ by $V$.  Let $\Gamma$ be the space of metrics
conformal to $g_0$ which have volume $V$.  Suppose that the operator
$A_{g_0}$ is {\it of type} $\Delta^{n/2}$, meaning that  it is  a
classical elliptic pseudodifferential operator on $M$ of degree $n$,
which is non-negative, self-adjoint with respect to $dV_0$, and has
null space precisely the constants, and moreover its leading order
term agrees with that of $\Delta^{n/2}$ where $\Delta$ is the
Laplace-Beltrami operator. For $g=F^{2/n}g_0$ in $\Gamma$, we set
$$
A_{ g}\ =\ F^{-1} A_{g_0}.
$$
Then trivially,  for any two metrics $g$ and $\tilde g=\tilde F^{2/n}g$ in $\Gamma$, we have
$$
A_{\tilde g}\ =\ \tilde F^{-1} A_{g}.
$$
For example, we can take  $A$ to be the GJMS operator $\square$
provided it is non-negative and has null space equal to the
constants, see [C], [FG], [GJMS], [GZ], [Gu]. However, we can work
with any operators $A_g$ constructed as above. Now fix $g\in\Gamma$
and let $dV$ denote the volume element for the metric $g$. The
operator $A=A_g$ is self-adjoint with respect to $dV$. Define
$A^{-1}$ to be the linear operator equal to the inverse of $A$ on
the orthogonal complement of the locally constant functions, and to
equal zero on the locally constant functions. The Green's function
for $A$  is the function on $M\times M$ which satisfies
$$
A^{-1}\phi\ =\ \int_M G(p,q) \phi(q)\,dV(q).
$$
Writing $d_g(p,q)$ for the Riemannian distance from $p$ to $q$ in the metric $g$, the function
$G(p,q)$ has an expansion at the diagonal of the form

$$
G(p,q)= -\frac{2n}{\gamma_n}\log d_g(p,q)\ +\ m(p)\ +\ O(d_g(p,q)),
$$
where
$$
\gamma_n\ =\ n!\ \omega_n,\qquad\qquad \omega_n=\text{ volume of
standard }n\text{ sphere }.
$$
The quantity $m(p)=m_g(p)$ is called the {\it Robin mass} at the point $p\in M$.
There are two
 natural ways to regularize the trace of $A^{-1}$.  One is to use the spectrum to take the
 {\it zeta regularization}.  It is simpler, however, to define
$$
\trace A^{-1}\ =\ \int_M m\,dV.
$$
For our results, the choice of definition is irrelevant since the two definitions differ by $c_n V$ where
$c_n$ is a universal constant depending only on dimension, see [Mor2], [St1], [St2], and also the
appendix of this paper.
  We note that recently Doyle and Steiner [DS1] gave a probabilistic interpretation
of $m(x)$ and $\trace \Delta^{-1}$ on closed surfaces.
\medskip

\proclaim{Theorem 1} Writing $\square_{S^n,V}$ for
the GJMS operator on the round $n$-sphere of volume $V$, we have
$$
\inf_{g\in \Gamma} \trace  A_g^{-1}\ \leq\ \trace \square_{S^n,V}^{-1}. \tag 1.4
$$
Moreover,  if $M$ is connected and the inequality (1.4) is strict,
then the infimum on the left hand side
 is attained on $\Gamma$ at a metric for which $m(p)$ is constant.
\endproclaim
\medskip

\noindent{\bf Remarks}.
1. Theorem 1 has some similarities with the Yamabe theorem, where the mass
is replaced by the scalar curvature, see for example [Sc], [Y]. In particular, we have
a lack of compactness in this problem which we get around using a method similar to [Y].

\noindent 2.  If the operator $A$ had trivial null space, then we would have
${
m_{F^{2/n}g} = m_g + \frac2{\gamma_n}\log F}$,
and the existence of constant mass metrics would be obvious.  This can
 be compared with  [H] where the existence of
constant mass metrics for the conformal Laplacian is immediate. The fact that
$A$ has non-trivial null space introduces a logarithmic HLS inequality into our analysis
of constant mass metrics.

\noindent 3. If the closed manifold $M$ has more than one component,
and if a minimizer for  $\trace \Delta^{-1}$ can be found for each
individual component when the volume of that component is fixed and
the metric varied  conformally,  then  a minimizer can be found on
$M$ by scaling the minimizers for each component appropriately. The
appropriate scaling can be found easily using a Lagrange multiplier.

\noindent 4. In [Ok3], it is shown that the inequality in (1.4) is
strict for any $2$-torus.
\medskip

From now on fix  a Riemannian metric $g$ on $M$ with volume element
$dV$, and for a positive function $F$ on $M$ define
$$
V_F\ =\ \int_M F\,dV.
$$
We introduce the space
$$
 (L\log L)^+(M)\ =\ \left\{ F:M\to [0,\infty)\ :\ F\text{ is measurable, }\int_M F\log F\,dV<\infty\right\}.
$$
\proclaim{Conformal change of $\trace A^{-1}$}  If $F$ is smooth and positive on $M$ then
$$
 \trace A_{F^{2/n}g}^{-1}\ =\ \mu(M,g,F), \tag 1.5
$$
where $\mu(M,g,F)$ is defined for $F\in (L\log L)^+(M)$ by
$$
 \mu(M,g,F)\ =\  \int_M  mF\,dV\ +\ \frac2{\gamma_n}\int_M F\log F\,dV
\ -\ \frac{1}{V_F}\int_M FA^{-1} F\,dV, \tag 1.6
$$
where $m$ is the mass for $A_g$.
\endproclaim
See [Mor1], [Mor2], [St1], [St2].  We give slightly different proof  of (1.5) in Section 3.  We can  rewrite
Theorem 1 in the following form.
\medskip

\proclaim{Theorem 1$'$} Let $M$ be a smooth $n$-dimensional manifold $M$
with Riemannian metric $g$.  Let $A$ be an operator on $M$ of type $\Delta^{n/2}$
with Robin mass $m$.    Then
$$
\inf\Sb V_F={V}\\ F\in (L\log L)^+(M)\endSb\mu(M,g,F)
\leq\
\trace\square_{S^n,V}^{-1}. \tag 1.7
$$
Moreover, if $M$ is connected and the inequality in (1.7) is strict,
the infimum on the left hand side is attained by some smooth
positive $F\in (L\log L)^+(M)$  with $V_F=V$ which satisfies the
pseudodifferential equation (or partial differential equation if $A$
is a partial differential operator)
$$
A(\log F)(p)\ =\ \gamma_n\left( \frac{F(p)}{V}\ -\ \frac{ Am(p)}2\ -\ 1\right).
\tag 1.8
$$
Furthermore, if $F=1$ attains the infimum on the left hand side of (1.7), then for
all $F\in (L\log L)^+(M)$ with $V_F=V$,
$$
\frac2{\gamma_n}\int_M F\log F\,dV
\ -\ \frac{1}{V}\int_M FA^{-1} F\,dV\ \geq\ 0. \tag 1.9
$$
\endproclaim
\medskip

\noindent{\bf Remarks}. 1. If $F$ is a local minimizer for $\mu(M,g,F)$ among functions
$F\in (L\log L)^+(M)$ with  $V_F=V$,
then $F$ is smooth and positive and satisfies (1.8). Just follow (2.20)-(2.22) with $\ep=0$.

\noindent 2. The functional $\mu(M,g,F)$ changes when either $g$ or $F$ is scaled.
We have stated (1.7) when $V_F=V$, but in view of (1.5) this is unnecessary.  Indeed, for any
positive constant $V'$ we have
$$
\inf\Sb V_F=V'\\ F\in (L\log L)^+(M)\endSb\mu(M,g,F)
\leq\
\trace\square_{S^n,V'}^{-1}.
$$
However,  we do need to assume $V_F=V$ in (1.9).
\medskip\medskip

\noindent{\bf Definition}. For $p\in M$, write
 $
 B_g(p,\delta)=\{q\in M:d_g(p,q)<\delta\}.
 $
 \medskip
  The inequality  (1.7) follows from   the next more precise result.
\proclaim{Theorem 2} Let $M$ be a smooth $n$-dimensional manifold $M$
with Riemannian metric $g$.  Let $A$ be an operator on $M$ of type $\Delta^{n/2}$
with Robin function $m$.  Then for each point $p\in M$,
\noindent
$$
\lim_{\delta\to 0}\ \inf\Sb
V_F=V\\ F\in (L\log L)^+(M)\\
\supp(F)\subset B_g(p,\delta) \endSb\mu(M,g,F)
\ =\
\trace\square_{S^n,V}^{-1}.
$$
Moreover, the limit is  uniform over $p\in M$.
\endproclaim
\medskip
The topology and local geometry of $M$ do not play any role in this
local result. The general  idea  is that if the function $F$ is
supported close to a point, then since locally all manifolds look
similar to $\Bbb R^n$, the degenerate metric $F^{2/n}g$ might as
well be on the sphere. Indeed,  $\mu(M,g,F)$ will be close to
$\trace \square^{-1}$ for some metric on the sphere. Moreover, by
choosing the function $F$ suitably we can ensure that
 the metric $F^{2/n}g$ blows up a neighborhood
of a point to be approximately
a round sphere, and  $\mu(M,g,F)$ will be close to $\trace \square_{S^n,V}^{-1}$. The details are
given in Section 2.
\medskip

\noindent{\bf Remark}.  In [DS2], a sequence of $2$ dimensional tori of fixed area
is constructed
in such a way that $\trace \Delta^{-1}$ converges to the value for the round sphere.
Whereas however, we make such a construction here for a fixed manifold by taking a sequence of conformal
factors which concentrate at a point, in [DS2] the sequence of tori is constructed by taking conformal
factors on a sequence of degenerating flat tori, and although the resulting tori approximate
the sphere, the conformal factors do not concentrate.
\medskip

Finally we show that the sharp form of the Logarithmic HLS inequality given in Theorem 1$'$ always
gives rise to  a sharp form of the Beckner-Onofri inequality. It is well known on $S^n$
that the logarithmic HLS inequalites and the Beckner inequalities  are dual to each other, see
[Be], [CL].
Those arguments can be extended to obtain a sharp Beckner-Onofri inequality in our situation.
 Suppose that $M$ is a smooth manifold with measure $d\s$ with $d\s(M)=1$, and
 $B:C^\infty(M)\to C^\infty(M)$
is a non-negative self-adjoint linear operator whose null space is equal to the constant functions, and
which is invertible on the orthogonal complement of the constants.

\proclaim{Theorem 3}  Let $\al$ be a smooth function on $M$ and let $\beta>0$ be a constant.
The following two statements are equivalent.

\noindent (a).  For $F\in C^\infty(M)$ with $F>0$ and $\int_M F\,d\s=1$,
$$
\frac\beta2\int_M F B^{-2} F\,d\s\ \leq\ \int_M F\log F\,d\s\ +\ \int_M F\al\,d\s.
$$

\noindent (b). For $u\in C^\infty(M)$,
$$
\frac1{2\beta}\int_M uB^2 u\,d\s\ \geq\
\log \int_M e^{u-\al}\,d\s\ -\ \int_M u\,d\s.
$$
\endproclaim

\proclaim{Corollary 4} (a). Suppose that the metric  $g\in \Gamma$ attains the infimum on the
left hand side of (1.7). Set $d\s=dV/V$.  Then for $u\in L^2_{n/2}(M)$,
$$
\frac{V}{2\gamma_n}\int_M uAu\,d\s\ \geq\
\log \left(\int_M e^{u}\,d\s\right)\ -\ \int_M u\,d\s.
$$
(b). Suppose that the inequality in (1.7) is an equality,
and $g$ is any metric in $\Gamma$. Set $d\s=dV/V$.  Then for  $u\in L^2_{n/2}(M)$,
$$
\frac{V}{2\gamma_n}\int_M uAu\,d\s\ \geq\
\log \left(\int_M e^{u-\gamma_n m/2}\,d\s\right)\ -\ \int_M u\,d\s\ +\ \frac{\gamma_n}{2V}
\trace \square_{S^n,V}^{-1}.
$$
\endproclaim
\medskip\medskip

\noindent{\bf Section 2.  Proofs of the Theorems}.
\medskip

\noindent{\bf Proof that Theorem 1 and Theorem $1'$ are equivalent.}
 In Theorem 1$'$, a metric
$g\in \Gamma$ is fixed, and other metrics in $\Gamma$ are expressed in the form
$F^{2/n}g$.  The value of $\trace A^{-1}$ is expressed in terms of $F$ by
 (1.5).  There are only three things to check.  Firstly,
 the fact that the metric $F^{2/n}g$ has constant Robin mass $m$ is equivalent to (1.8).
 This follows from Lemma 2.1 below which is proved in [St1], [St2].  We give another proof in Section 3.

\proclaim{Lemma 2.1}
{\bf Conformal change of the Robin mass}.
 Suppose $g$ is a metric on $M$ and $F$ is a smooth positive function on $M$. Write $m_g(p)$
 for the Robin constant of $A_g$ at $p$.  Then
$$
m_{F^{2/n}g}(p)\ =\  m _g(p)\ +\ \frac{2\log F}{\gamma_n}\ -\  \frac{2}{ V_F} (A^{-1}F)(p)
\ +\  \frac1{V_F^2}\int_M F
(A^{-1}F)\,dV. \tag 2.1
$$
\endproclaim

The second point  that needs to be clarified is that in (1.4) we
take an infimum over the set of smooth metrics $\Gamma$, while  in
(1.7) we allow degenerate metrics of the form $F^{2/n}g$ where the
non-negative function $F$ is in $(L\log L)^+(M)$. The fact that this
does not change the infimum is due to the fact that the smooth
positive functions are dense in $(L\log L)^+(M)$ and the functional
$F\to \mu(M,g,F)$ is continuous on $(L\log L)^+(M)$.  To see the
latter, we need to show that $F\to \int_M FA^{-1} F\,dV$ is
continuous, which follows easily from   the next Lemma.

\proclaim{Lemma 2.2} {\bf A simple logarithmic Sobolev inequality.}

\noindent For each $\ep>0$, there exists a constant
$C_\ep=C_\ep(M,g)$ such that for all $F\in (L\log L)^+(M)$,
$$
   \|A^{-1}F(x)\|_\infty \leq\ (1+\ep)\frac2{\gamma_n}\int_M F\log F\,dV\ +\ C_\ep \left(\int_M F\,dV
   +1\right). \tag2.2
$$
\endproclaim
\noindent{\bf Remarks}. In fact $A^{-1}F$ is continuous when  $F\in (L\log L)^+(M)$.
There is no bound of the form (2.2)  when $\ep=0$.
\medskip

This completes the proof that Theorem 1 and Theorem 1$'$
are equivalent, except we should just note that although (1.9) does not feature in Theorem 1, it
 follows easily from (1.8).  Indeed, if $F=1$ attains the minimum
of the left hand side of (1.7), then from (1.8) we see that $m=m(x)$ must be constant.  Then
 for $F\in (L\log L)^+(M)$ with $V_F=V$,
$$
mV\ =\ \mu(M,g,1)\ \leq\ \mu(M,g,F)\ =\ mV\ +\ \frac2{\gamma_n}\int_M F\log F\,dV
\ -\ \frac{1}{V_F}\int_M FA^{-1} F\,dV.
$$
\medskip

\noindent{\bf Proof of Theorem 2.}  It will be convenient to give a characterization of
$\trace \square_{S^n,V}^{-1}$ involving the flat space $\Bbb R^n$ rather than the sphere.
This equivalence can be proved by using stereographic projection to identify $\Bbb R^n$ with
$S^n$ as in [CL].    We give  a natural proof here using  concentration arguments,
see Proposition 2.6.
Following previous notation, for a positive function $f$ on $\Bbb R^n$ we set
$$
V_f\ =\ \int_{\Bbb R^n} f\,dx.
$$
Define
$$
(L\log L)^+_c(\Bbb R^n)\ =\ \left\{ f:\Bbb R^n\to [0,\infty)\ :\ f\text{ is measurable with compact
support, }
\int_{\Bbb R^n} f\log f\,dx<\infty\right\}.
$$
and for $f\in (L\log L)_c^+(\Bbb R^n)$  define the functional on $\Bbb R^n$ analogous to $\mu(M,g,F)$:
$$
\mu(\Bbb R^n,f)\ =\
\frac2{\gamma_n}\left( \int_{\Bbb R^n} f\log f\,dx\ -\ \frac{n}{V_f}\int_{\Bbb R^{2n}}
f(x)\log |x-y| f(y)\,dxdy\right).
$$
Set
$$
\mu_V(\Bbb R^n)\ =\ \inf\Sb V_f=V\\ f\in (L\log L)^+_c(\Bbb R^n)\endSb\mu(\Bbb R^n,f).
$$
\medskip

We will now prove Theorem 2 with the quantity $\trace \square_{S^n,V}^{-1}$ replaced by $\mu_V(\Bbb R^n)$.
It will then emerge that these quantities are equivalent.
\proclaim{Lemma 2.3}
For each point $p\in M$,
\noindent
$$
\lim_{\delta\to 0}\ \inf\Sb
V_F=V\\ F\in (L\log L)^+(M)\\
\supp(F)\subset B_g(p,\delta) \endSb\mu(M,g,F)
\ =\
\mu_V(\Bbb R^n).
$$
and the limit is uniform over $p\in M$.
\endproclaim

Lemma 2.3 follows immediately from the next Lemma.

\proclaim{Lemma 2.4}  Given $\ep>0$, there exists $\delta>0$ such that:

\noindent (a). If $F$ is supported in
$B_g(p,\delta)$ for some $p\in M$, and $F\in (L\log L)^+(M)$,  then there exists
$f\in (L\log L)^+_c(\Bbb R^n)$ with $V_f=V_F$ such that
$$
|\mu(M,g,F)\ -\ \mu(\Bbb R^n,f)|\ \leq\ \ep V_F.
$$

\noindent (b). If $f\in (L\log L)^+_c(\Bbb R^n)$
and $p\in M$,  there exists   $F\in (L\log L)^+(M)$ supported in $B_g(p,\delta)$
with $V_F=V_f$ such that
$$
|\mu(M,g,F)\ -\ \mu(\Bbb R^n,f)|\ \leq\ \ep V_f.
$$
\endproclaim

Notice that  from (a) we get
$$
\lim_{\delta\to 0}\inf\Sb
V_F=V\\ F\in (L\log L)^+(M)\\
\supp(F)\subset B_g(p,\delta) \endSb\mu(M,g,F) \ \geq\  \mu_V(\Bbb
R^n),
$$
and from (b) we get
$$
\lim_{\delta\to 0}\inf\Sb
V_F=V\\ F\in (L\log L)^+(M)\\
\supp(F)\subset B_g(p,\delta) \endSb\mu(M,g,F) \ \leq\ \mu_V(\Bbb
R^n).
$$

\noindent{\bf Proof of Lemma 2.4}. The idea of the proof is that
if $F\in (L\log L)^+(M)$ is supported close to a point $p_0\in M$, by taking
suitable coordinates we can  identify
$F$ with a function $f$ on $\Bbb R^n$ so that    $\mu(M,g,F)$  is close to $\mu(\Bbb R^n,f)$.
Conversely, given $f\in (L \log L)_c^+(\Bbb R^n)$ we can rescale so that the support of $f$  becomes
small while $V_f$ and $\mu(\Bbb R^n,f)$ remain constant, and then we can consider the function
$F$ on $M$ which is given in coordinates by $f$.
To carry out the details, for arbitrary fixed $\ep>0$,  choose $\delta>0$ so
that if $d_g(p,q)<2\delta$ then
$$
\left|G(p,q)\ +\ \frac{2n}{\gamma_n}\log d_g(p,q)\ -\ m(p)\right|\ <\ \ep.\tag 2.3
$$
We now choose good coordinates on $M$ around each point $p_0$.

\proclaim{Lemma 2.5}  For each $\ep>0$, there exists $\delta>0$ such that for each $p_0\in M$ there
exists an open neighborhood $U$ of $0$ in $\Bbb R^n$, and
smooth coordinates $x=(x_1,\dots,x_n):B_g(p_0,\delta)\to U$ such that
 $x(p_0)=0$,
 $$
 dV=dx=dx_1\dots dx_n, \tag 2.4
 $$
and if $|x|$ is the Euclidean norm
 of $x\in \Bbb R^n$, then for $p,q\in B_g(p_0,\delta)$,
 $$
 e^{-\ep} |x(p)-x(q)|\ \leq\ d_g(p,q)\ \leq\ e^\ep|x(p)-x(q)|.\tag 2.5
 $$
 \endproclaim
 This lemma is proved in Section 3.

Now given $\ep>0$  choose $\delta$  small enough so that (2.3) and the conclusion of Lemma 2.5 hold,
and for $p_0\in M$ take the coordinates $(x_1,\dots,x_n)$  of Lemma 2.5.
 Then if $p,q\in B_g(p_0,\delta)$, we have
$$
G(p,q)\ =\ -\frac{2n}{\gamma_n}\log |x(p)-x(q)|\ +\ m(x(p))\ +\ \eta(p,q),
\qquad\quad |\eta(p,q)|<c\ep,\qquad\quad c=1+\frac{2n}{\gamma_n}.
$$
If $F\in (L\log L)^+(M)$ is supported  in $B(p_0,\delta)$, then
$$
\multline
\frac1{V_F} \int_M FA^{-1}F\,dV\ =\ \frac1{V_F}\int_M\int_M F(p)G(p,q)F(q)\,dV(p) dV(q)
\\ =\ -\frac{2n}{\gamma_n V_F}\int_M\int_M F(p)\log|x(p)-x(q)|F(q)\,dV(p) dV(q)
\ +\ \int_M mF\,dV\\ +\ \frac1{V_F} \int_M\int_M F(p)\eta(p,q) F(q)\,dV(p) dV(q).
\endmultline
\tag 2.6
$$
Define the function $f$ on $U$ by $f(x(p))=F(p)$, and extend
$f$ to $\Bbb R^n$ by setting it equal to zero outside $U$. Then from (2.6),
$$
\multline
\mu(M,g,F)\ =\  \int_M  mF\,dV\ +\ \frac2{\gamma_n}\int_M F\log F\,dV
\ -\ \frac{1}{V_F}\int_M FA^{-1} F\,dV\\ =\ \mu(\Bbb R^n,f)\  +\
\frac1{V_F} \int_M\int_M F(p)\eta(p,q) F(q)\,dV(p) dV(q).
\endmultline
\tag 2.7
$$
The second term on the right is bounded by $c\ep V_F$. If we replace $\ep$ by $\ep/c$
we get (a).
To prove (b), we note that
the functional $\mu(\Bbb R^n,f)$ has a scale invariance.  Indeed,
for $f\in (L\log L)^+_c(\Bbb R^n)$ and  $\la>0$, we can set
$$
h(x)\ =\ \frac1{\la^n}\  f\left(\frac{x}\la\right) \tag 2.8
$$
and then one can check that  $V_h=V_f$ and $\mu(\Bbb R^n,h)=\mu(\Bbb R^n,f)$.
Now given $\ep>0$, we pick $\delta>0$ so that (2.3) and the conclusion of Lemma 2.5 hold,
and given $p_0\in M$, choose the coordinates $(x_1,\dots,x_n)$ on $B_g(p_0,\delta)$
from Lemma 2.5. Then the image of $B_g(p_0,\delta)$ under the coordinates map contains some ball
$B_{\Bbb R^n}(0,\delta')=\{y\in \Bbb R^n:|y|<\delta'\}$.  For each function $f\in (L\log L)^+_c(\Bbb R^n)$
 we can choose $\la$ sufficiently small so that $h(x)$ defined by (2.8) is supported in
 $B_{\Bbb R^n}(0,\delta')$.
 But then define the function $F$ on $M$ supported in $B_g(p_0,\delta)$ by
 $F(p)=h(x(p))$.  Then (2.7) holds as before and we get (b).

To prove Theorem 2, it remains to show
\proclaim{Proposition 2.6}
$$
 \mu_V(\Bbb R^n)\ =\ \trace \square_{S^n,V}^{-1}.
$$
\endproclaim
To prove this
we first notice that by applying Lemma 2.3 to the sphere $S^n$ with the round metric $g_0$
of volume $V$, we have
$$
\mu_V(\Bbb R^n)\ =\ \lim_{\delta\to 0}\ \inf\Sb
V_F=V\\ F\in (L\log L)^+(S^n)\\
\supp(F)\subset B_g(p,\delta) \endSb\mu(S^n,g_0,F) \tag 2.9
$$
On the other hand, by Morpurgo's interpretation of the logarithmic HLS inequality,
$$
\trace \square_{S^n,V}^{-1}\ =\  \inf\Sb
V_F=V\\ F\in (L\log L)^+(S^n) \endSb\mu(S^n,g_0,F). \tag 2.10
$$
Since the infimum in (2.10) is over a bigger class of functions than in (2.9), we see that
$$
\trace \square_{S^n,V}^{-1}\ \leq\ \mu_V(\Bbb R^n). \tag 2.11
$$
The equality in (2.10) is obtained exactly when $F$ is the Jacobian
of a conformal transformation of the sphere.  Now naturally there is
no sequence of such functions $F_j$ whose supports shrink down to a
point.  However, it is a well known that by combining stereographic
projection with  the  scaling of Euclidean space, there exists a
sequence of conformal transformations of $S^n$ whose Jacobians $F_j$
concentrate in the following sense.
\medskip

\noindent{\bf Definition}.  Suppose that $F_j$ is a sequence of
positive functions in $L^1(M))$, we say that {\it $F_j$ concentrates
 as $j\to \infty$}, if for every $\delta>0$, there exists $j_0$
such that if $j\geq j_0$ then there exists $p_j\in M$ with
$$
\int_{B_g(p_j,\delta)} F_j\,dV\ >\ V_{F_j}(1-\delta).
$$
\medskip
 Proposition 2.6 follows from the next result which is an extension of Proposition 2.3.
This then completes the proof of Theorem 2.

\proclaim{Proposition 2.7}  Let $g_0$ be the round metric of volume $V$ on $S^n$.
If $F_j\in (L\log L)^+(S^n)$ is a concentrating sequence with $V_{F_j}=V$, then
$$
\liminf_{j\to\infty} \mu(S^n,g,F_j)\ \geq\ \mu_V(\Bbb R^n).
$$
\endproclaim
Proposition 2.7 is proved in Section 3. In fact later on we will need the same result for general manifolds
$M$, see Proposition 2.12.
\medskip\medskip

\noindent{\bf Proof of Theorem 1$'$}.  To prove Theorem 1$'$ we will need
to start with a weak form given in Proposition 2.9, and we will prove this
by using the fact that it is true on $S^n$ (or equivalently the fact that $\mu_V(\Bbb R^n)$
is finite) and moving this inequality over to $M$ using a partition of unity.
To accomplish this we need the following polarized form.
\medskip

\proclaim{Lemma 2.8} Let $g$ be a metric on $M$ with volume $V$, and let
$A$ be an operator on $M$ of type $\Delta^{n/2}$. For a function $F$ on $M$, set
$\s_F=V_F/V$.
The following three statements are equivalent.

\noindent (a).  For every $F\in (L\log L)^+(M)$ with $V_F=V$,
$$
\frac1V\int_M F A^{-1} F\,dV\ \leq\ \frac2{\gamma_n}\int_M F \log F\,dV\ +\ C.
$$

\noindent (b).  For all $Q,R\in (L\log L)^+(\Bbb R^n)$
with $V_Q=V_R=V$,
$$
\frac1V\int_M Q A^{-1} R\,dV\ \leq\ \frac1{\gamma_n}\left(\int_M Q\log Q\,dV + \int_M R\log R\,dV\right)
\ +\ C.
$$

\noindent (c). For all $Q,R\in (L\log L)^+(\Bbb R^n)$,
$$
\frac1V\int_M Q A^{-1} R\,dV\ \leq\ \frac1{\gamma_n}\left( \s_R\int_M Q\log \frac{Q}{\s_Q}\,dV
\ +\  \s_Q \int_M R\log \frac{R}{\s_R}\,dV\right)\ +\ C\,\s_Q\,\s_R.
$$
\endproclaim
\medskip

\proclaim{Proposition 2.9}  ({\bf Weak logarithmic HLS inequality}.)
There exists a constant $C=C(M,g)$ such that if $F\in (L\log L)^+(M)$ with $V_F=V$, then
$$
 \frac{2}{\gamma_n}\int_M F\log F\,dV
 \ -\  \frac1{V}\int_M F A^{-1} F\,dV\ \geq\  -C.\tag 2.12
$$
\endproclaim
Lemma 2.8 and Proposition 2.9 are proved in Section 3.
In order to prove Theorem 1$'$, we will need to get around a lack of compactness, which we
do by adapting the ideas of [Y] to our situation.  Let $\la_1$ be the lowest positive eigenvalue
of $A$, and consider the functional
$$
\mu^{(\ep)}(M,g,F)\ =\  \int_M  mF\,dV\ +\ \frac2{\gamma_n}\int_M F\log F\,dV
\ -\ \frac{(1-\ep)\la_1^\ep}{V_F}\int_M FA^{-1-\ep} F\,dV. \tag 2.13
$$
\proclaim{Proposition 2.10} There exist a positive function $F^{(\ep)}\in C^\infty(M)$
with $V_{F^{(\ep)}}=V$ which minimizes $\mu^{(\ep)}$, that is
$$
\mu^{(\ep)}(M,g,F^{(\ep)})\ =\ \inf\Sb V_F=V\\ F\in (L\log L)^+(M)\endSb \mu^{(\ep)}(M,g,F). \tag 2.14
$$
Moreover, $F^{(\ep)}$ satisfies the equation
$$
A^{1+\ep}(\log F^{(\ep)})\ =\ \gamma_n\left( \frac{(1-\ep)\la_1^\ep F^{(\ep)}(x)}V
\ -\ \frac{ A^{1+\ep}m(x)}2\ -\ (1-\ep)\la_1^\ep\right).
\tag 2.15
$$
\endproclaim

\noindent {\bf Proof of Proposition 2.10}.  Applying Proposition 2.9, we see that if $V_F=V$, then
$$
\multline
\mu^{(\ep)}(M,g,F)\ \geq\ \int_M  mF\,dV\ +\ \frac2{\gamma_n}\int_M F\log F\,dV
\ -\ \frac{1-\ep}{V}\int_M FA^{-1} F\,dV\\ =\ \frac{2\ep}{\gamma_n}\int_M
F\log F\,dV\ +\ (1-\ep)\left(\frac{2}{\gamma_n}\int_M F\log F\,dV
 \ -\  \frac1{V}\int_M F A^{-1} F\,dV\right)
\ +\  \int_M mF\,dV\\ \geq \ \frac{2\ep}{\gamma_n}\int_M
F\log F\,dV\ -\ C'
\endmultline
\tag 2.16
$$
Choose a sequence $F_j$  in $(L\log L)^+(M)$ with $V_{F_j}=V$ for all $j$ and
$$
\lim_{j\to\infty} \mu^{(\ep)}(M,g,F_j)\ =\ \inf\Sb V_F=V\\ F\in (L\log L)^+(M)\endSb \mu^{(\ep)}(M,g,F).
\tag 2.17
$$
Then by (2.16), there exists $C$ independent of $j$ such that
$$
\int_M F_j \log F_j\,dV\ \leq\ C.
$$

\proclaim{Lemma 2.11}  Suppose that $\Phi:[0,\infty)\to \Bbb R$ is a continuous convex function
with
$$
 \frac{\Phi(t)}t\ \to \ \infty\qquad \text{ as }t\to\infty.
$$
Suppose that $F_j$ is a sequence of non-negative measurable functions on $M$ such that
$$
\sup_j \int_M \Phi(F_j)\,dV\ =\ S\  <\ \infty.
$$
Then after replacing $F_j$ by a subsequence, there exists $F\in L^1(M)$ such that
$F_j\to F$ weakly, that is for every $\phi\in C(M)$,
$$
\int_M F_j \phi\,dV\ \to\ \int_M F\phi\,dV\qquad \text{ as }j\to \infty,
$$
and
$$
\int_M \Phi(F)\,dV\ \leq \ \liminf_{j\to \infty} \int_M \Phi(F_j)\,dV. \tag 2.18
$$
\endproclaim
We will prove this Lemma in Section 3. Applying it
 to the function $\Phi(t)=t\log t$ and the sequence $F_j$ in (2.17),
 we see that by taking a subsequence of $F_j$ if necessary
there exists
$F^{(\ep)}\in (L\log L)^+(M)$
with $F_j\to F^{(\ep)}$ weakly in $L^1(M)$ and
$$
\int_M F^{(\ep)}\log F^{(\ep)}\,dV\ \leq\ \liminf_{j\to \infty} F_j\log F_j\,dV.
\tag 2.19
$$
However, the integral kernel of $A^{-1-\ep}$ is continuous on $M\times M$, and so by simple estimates
$A^{-1-\ep}$ is
bounded from $L^1(M)$ to $C(M)$, and furthermore $\{A^{-1-\ep} F:\|F\|_1\leq C\}$ is equicontinuous.  Hence
by taking a subsequence we can assume
 $A^{-1-\ep}F_j$ converges in  $C(M)$.  But since $F_j$ converges weakly in $L^1(M)$ to $F^{(\ep)}$, this
 limit  must equal $A^{-1-\ep}F^{(\ep)}$.  Then
$$
\lim_{j\to\infty } \int_M F_j A^{-1-\ep} F_j\,dV\ =\ \int_M F^{(\ep)}A^{-1-\ep} F^{(\ep)}\,dV.
$$
But from this and (2.19) we get
$$
\mu^{(\ep)}(M,g,F^{(\ep)})\ \leq\ \lim_{j\to\infty} \mu^{(\ep)}(M,g,F_j),
$$
and hence  $\mu(M,g,F^{(\ep)})$ satisfies (2.14). Now we need to
show that  $F^{(\ep)}$ is positive and smooth. To simplify notation,
write $F=F^{(\ep)}$ and $B=(1-\ep)\la_1^\ep A^{-1-\ep}$. We will
first show that $F$ is bounded below by a positive constant. For
every bounded function $H$ such that $F+H\in (L\log L)^+(M)$ and
$\int_M H\,dV=0$, we have
$$
\align
&\mu^{(\ep)}(M,g,F+H)-\mu^{(\ep)}(M,g,F)\ =\ \int_M mH\,dV
\ +\ \frac2{\gamma_n}\int_M \left((F+H)\log (F+H)-F\log F\right)\,dV
 \\ &\qquad\qquad\qquad\qquad\qquad\qquad\qquad\quad
 -\ \frac{2}{V}\int_M H B F\,dV\ -\ \frac{1}V \int_M H B H\,dV \tag 2.20\\
 \allowdisplaybreak &\leq\ \int_M mH\,dV
\ +\ \frac2{\gamma_n}\int_M \left((F+H)\log (F+H)-F\log F\right)\,dV
\ -\ \frac{2}{V}\int_M H B F\,dV\\
 \allowdisplaybreak & \leq\ \frac2{\gamma_n}\int_M \left((F+H)\log (F+H)-F\log F\right)\,dV
 \ +\ C\int_M|H|\,dV,\tag 2.21
 \endalign
$$
where
$$
C\ =\ \|m\|_\infty\ +\ \frac2V \|BF\|_\infty .
$$
We will show that if $F$ is very small on a set of positive measure, then by choosing
$H$ appropriately, we can make (2.20) negative
hence contradicting the fact that $F$ minimizes $\mu^{(\ep)}$. By the mean value theorem,
$$
(f+h)\log (f+h)\ -\ f\log f\ <\ h(1+\log (f+h))\ <\  0
$$
whenever
$$
f>0,\, h>0,\,\quad f+h<1/e,\qquad\qquad\text{or}\qquad\qquad     f>0,\, h<0, \, f+h> 1/e.
$$
Since the mean value of $F$ equals $1$, the set  $U=\{p\in M:F(p)>1/e\}$ has positive
measure.   Suppose that $N>1$ and the set $U_N=\{ p\in M:F(p)<e^{-N}\}$ has positive measure.  Define
$H$ so that $\int_M H\,dV=0$ and
$$
\cases 0<H<e^{-N}-F \qquad & \text{ on }U_N \\ e^{-1}-F<H<0 \qquad &\text{ on }U\\
H=0 & \text{ on } M\setminus (U_N\cup U).\endcases
$$
Then by (2.21),
$$
\int_M \left( (F+H)\log (F+H)-F\log F\right)\,dV\ \leq\ -N\int_{U_N}H\,dV\ =\ \frac{-N}2\int_M|H|\,dV.
$$
From (2.21), we see that if $N/\gamma_n >C$ then $\mu^{(\ep)}(M,g,F+H)-\mu^{(\ep)}(M,g,F)<0$, which
contradicts $F$ being a local minimum.  Hence $F$ is bounded below.  A similar argument shows that
$F$ is bounded above. Setting the variation of  (2.20) about
$H=0$ equal to zero, shows that
$$
m\ +\ \frac{2\log F}{\gamma_n}\ -\ \frac{2B F}{V}\ =\ \text{constant}. \tag 2.22
$$
However, since $\log F$ is in $L^\infty(M)$ and $m$ is smooth, we can apply elliptic regularity
theory to conclude that $F$ is smooth.  Applying the operator $A^{-1-\ep}$ to (2.22) gives (2.15) and
completes the proof of Proposition 2.10.
\medskip\medskip

We return to the proof of Theorem 1$'$. We want to obtain a
minimizer for $\mu(M,g, F)$. Take a sequence $\ep_j\searrow 0$, and
consider the sequence of functions $F^{(\ep_j)}$. Since
$F^{(\ep_j)}$ is minimizing for $\mu^{(\ep)}(M,g,F)$ and
$\mu^{(\ep_j)}(M,g,F)$ decreases to $\mu(M,g,F)$ as $\ep_j\searrow
0$, we see that
$$
\lim_{j\to\infty }\mu(M,g,F^{(\ep_j)})\ =\ \inf\Sb V_F=V\\ F\in (L\log L)^+(M)\endSb \mu(M,g,F).\tag 2.23
$$
\medskip
We will need the analog of Proposition 2.7 on the  manifold $M$.

\proclaim{Proposition 2.12}
If $F_j\in (L\log L)^+(M)$ is a concentrating sequence with $V_{F_j}=V$, then
$$
\liminf_{j\to\infty} \mu(M,g,F_j)\ \geq\ \mu_V(\Bbb R^n).
\tag 2.24
$$
\endproclaim
\medskip
Proposition 2.12 will be proved in Section 3.  Applying Propositions 2.12 and 2.6 to the sequence
$F^{(\ep_j)}$, we see that if $F^{(\ep_j)}$ concentrates, then
$$
\inf\Sb V_F=V\\ F\in (L\log L)^+(M)\endSb \mu(M,g,F)\ =\ \mu_V(\Bbb R^n)\ =\ \trace \square_{S^n,V}^{-1}.
$$

On the other hand if the sequence $F^{(\ep_j)}$ does not concentrate
then after taking a subsequence,
 there exists  $\delta>0$ such that for every $p\in M$ and every $j$,
$$
\int_{B(p,\delta)} F^{(\ep_j)}\,dV\ \leq\ V(1-\delta).
$$
We can apply the following result which is proved in Section 3.

\proclaim{Proposition 2.13}  ({\bf Improved logarithmic HLS inequality for non-concentrating functions}.)
 Given $\delta>0$, there exists a constant $C=C(M,g,\delta)$ such that
if   $F\in (L\log L)^+(M)$ satisfies $V_F=V$ and is such that
for every $p\in M$,
$$
\int_{B_g(p,\delta)} F\,dV\ \leq\ (1-\delta)V_F,
$$
then
$$
 (1-\delta)\frac{2}{\gamma_n}\int_M F\log F\,dV
 \ -\  \frac1{V_F}\int_M F A^{-1} F\,dV\ \geq\  -C.
$$
\endproclaim
From this result and the fact that $\mu(M,g,F^{(\ep_j)})$ is bounded above, we see that there is a uniform bound
$$
\int_M F^{(\ep_j)}\log F^{(\ep_j)}\,dV\ \leq\ C, \tag 2.25
 $$
\proclaim{Lemma 2.14}  The operators $A^{-\ep}$ with $\ep\in [0,1/2]$
are uniformly bounded on $C^k(M)$.
\endproclaim
We sketch the proof of Lemma 2.14 in Section 3.
Applying  Lemmas 2.2 and 2.14 to (2.25),  we get a constant $C$ such that
$$
 \|A^{-1-\ep_j}F^{(\ep_j)}\|_\infty\ \leq\ C
$$
From (2.15), this gives a uniform bound on $\|\log F^{(\ep)}\|_\infty$.  Indeed,
since the null space of $A$ is the constant functions, we  get  a bound of the form
$$
\| \log F^{(\ep)}-c_\ep\|_\infty\ \leq\ C',
$$
where $C'$ is a constant independent of $\ep$ but $c_\ep$ is an unknown constant which may depend
on $\ep$.  However, since the average value of $F^{(\ep)}$ is $1$, we see that $F^{(\ep)}$
takes the value $1$ and so $\log F^{(\ep)}$ takes the value zero.  Hence $|c_\ep|\leq C'$.
Hence we get uniform bounds  on $\|F^{(\ep_j)}\|_\infty$.
However, $A^{-1}$ is bounded from $C^k$ to $C^{k+n-1}$ and so applying (2.15) again, we obtain uniform
bounds on $\|F^{(\ep_j)}\|_{C^{n-1}(M)}$.  Continuing in this way, for all $k$ we get  bounds on
$\|F^{(\ep_j)}\|_{C^k(M)}$ which are uniform in $j$. Hence using the Azela-Ascoli theorem and
and a diagonalization argument, we can find a subsequence
which converges in $C^\infty(M)$ to a smooth positive function $F$ which attains the right hand
side of (2.23), and which  satisfies the limiting equation (1.8).
\medskip\medskip

\noindent{\bf Proof of Theorem 3}.   Set $\s_u=\int_M u\,d\s$. Our starting point is the inequality
$$
\align
 \int_M F(u-\s_u)\,d\s\ &\leq\ \left(\int_M |B^{-1} F|^2\,d\s\right)^{1/2}
\left(\int_M |B^{} u|^2\,d\s\right)^{1/2}\\ &\leq\ \frac{\beta}2 \int_M FB^{-2}F\,d\s\ +\
\frac1{2\beta}\int_M uB^2 u\,d\s, \tag 2.26
\endalign
$$
with equality if $F=c+B^2u/\beta$, where $c$ is a constant.

\noindent (a)$\Rightarrow$(b). We follow [CL].
Assume that (a) holds. Then for $u\in C^\infty(M)$,
applying (2.26) with  $F> 0$ and $\int_M F\,d\s=1$, we have
$$
\align
\frac1{2\beta}\int_M uB^2 u\,d\s\ +\ \int_M u\,d\s\ &\geq \  \int_M Fu\,d\s\ -\ \frac{\beta}2 \int_M FB^{-2} F\,d\s\\
&\geq\ \int_M F(u-\al)\,d\s\ -\ \int_M F\log F\,d\s.
\endalign
$$
We get (b) by choosing
$$
F\ =\ \frac{e^{u-\al}}{\int_M e^{u-\al}\,d\s}.
$$

\noindent (b)$\Rightarrow$(a).  Assume that (b) holds, and suppose $F\in C^\infty(M)$ satisfies
$F> 0$ and $\int_M F\,d\s=1$, and set $u=\beta B^{-2}F$.  Then
$$
\align
\frac{\beta }2 \int_M FB^{-2}F\,d\s\ &=\ \int_M F u\,d\s\ -\ \frac1{2\beta} \int_M uB^2u\,d\s\\
&\leq\ \int_M Fu\,d\s\ -\ \log \int_M e^{u-\al}\,d\s\\
&\leq\ \int_M F\log F\,d\s\ +\ \int_M F\al\,d\s,
\endalign
$$
where the last line follows from Jensen's inequality:
$$
\int_M (u-\al-\log F)F\,d\s\ \leq\ \log \int_M e^{u-\al-\log F} F\,d\s\ =\ \log \int_M e^{u-\al}\,d\s.
$$

\noindent{\bf Remark}. In proving (a)$\Rightarrow$(b), the choice $F=e^{u-\al}/\int_M e^{u-\al}\,d\s$
is the Legendre
function for the functional $\int_M F\log F\,d\s$.
\medskip\medskip

\noindent{\bf Section 3.  Proofs of the auxiliary results.}
\medskip

\noindent{\bf Proof of Lemma 2.1}.  Write $d(p,q)$ for the distance from $p$ to $q$ in the
metric $g$.
The Green's function $G(p,q)$ for $A=A_g$ is a smooth function
of $(p,q)\in M\times M$ away from the diagonal $p,q$, and is characterized by
the following conditions:
$$
\align
& {A}_q G(p,q)=-\frac1V\qquad \text{ for }q\neq p, \\
& G(p,q)\ +\ \frac{2n}{\gamma_n} \log d(p,q)\qquad\text{ is bounded},\tag 3.1\\
&\int_M G(p,q)\,dV(q)\ =\ 0,\tag 3.2\\
& G(p,q)=G(q,p).
\endalign
$$
Here, $A_q G(p,q)$ denotes the operator $A$ applied to $G(p,q)$ in the variable $q$.
Consider the function
$$
E(p,r,q)=G(p,q)-G(r,q)
$$
which satisfies
$$
\align
& {A}_q E(p,r,q)\ =\ 0,\qquad\text{ for } q\neq p,r,  \tag 3.3\\
& E(p,r,q)\ +\ \frac{2n}{\gamma_n}\left( \log d(p,q)-\log d(r,q)\right)\quad \text{is bounded}, \\
&\int_M E(p,r,q)\,dV(q)=0.
\endalign
$$
Now set $\tilde g=F^{2/n}g$. We write $\tilde G$ for the  Green's function for $\tilde g$, and
$$
\tilde E(p,r,q)=\tilde G(p,q)-\tilde G(r,q).
$$
 Because of (3.3) we see that
$$
{A}_q(\tilde E(p,r,q)-E(p,r,q))=0.
$$
By (3.1) we see that
$\tilde E(p,r,q)-E(p,r,q)$ is bounded on $M$.  We conclude  that $\tilde E(p,r,q)-E(p,r,q)$ is constant.
We can compute this constant by applying (3.2) for $\tilde E$, and we get
$$
\tilde G(p,q)-\tilde G(r,q)\ =\ G(p,q)-G(r,q)\ -\ \frac1{V_F} (A^{-1}
F)(p) \ +\ \frac1{V_F}(A^{-1} F)(r).
$$
 Averaging  with respect to $F(r) dV(r)$, we see
that
$$
\tilde G(p,q)\ =\ G(p,q)\   -\ \frac1{V_F} (A^{-1}F)(q) \ -\
\frac1{V_F} (A^{-1}F)(p) \ +\  \frac1{V_F^2}\int_M
F A^{-1}F\,dV
$$
Now
$$
m(p)\ =\ \lim_{q\to p}\left(  G(p,q)+\frac{2n}{\gamma_n}\log d(p,q)\right),
$$
writing $\tilde d$ for the distance function for the metric $\tilde g$,
we see from (3.12) that
$$
\multline  m_{F^{2/n}g}(p)\ =\ \lim_{q\to p}\biggl(\frac{2n}{\gamma_n}\log \tilde d(p,q)\ -\
 \frac{2n}{\gamma_n}\log d(p,q)\ +\  m_g (p) \\ -\
\frac{1}{ V_F} (A^{-1}F)(q)
\ -\ \frac{1}{ V_F}  (A^{-1}F)(p)
\ +\  \frac{1}{ V_F^2}\int_M F A^{-1}F\,dV,\biggr)
\endmultline
$$
so
$$
 m_{F^{2/n}g}(p) \ =\  m_g(p) \ +\ \frac{2\log F(p)}{\gamma_n}\ -\  \frac{2}{ V_F} A^{-1}F (p)
\ +\  \frac1{ V_F^2}\int_M F
A^{-1}F\,dV.
$$
Hence defining
$$
\trace A_{ g}^{-1}\ =\ \int_M m\, dV,
$$
we have
$$
\trace A_{F^{2/n}g}^{-1}\ =\ \int_M m F\,dV\ +\ \frac2{\gamma_n}\int_M F\log F\,dV\
-\
  \frac{1}{ V_F}\int_M F (A^{-1}F)\,dV.
$$
\medskip\medskip

\noindent{\bf Proof of Lemma 2.2}.  Writing $G(p,q)$ for the Green's function for $A_g$, we have
$$
G(p,q)\ =\ -\frac{2n}{\gamma_n}\log d(p,q)\ +\ E(p,q),
$$
where $E(p,q)$ is bounded.  Hence
$$
A^{-1} F(p)\ =\ \int_M G(p,q)F(q)\,dV(q)\ =\ \frac{-2n}{\gamma_n}\int_M F(q) \log d(p,q) \,dV(y)\ +\
\int_M E(p,q) F(q)\,dV(q),
$$
and the second term on the right is bounded by
$V\sup E$, which can be absorbed in $C_\ep$.  For the first term,  taking $0<\delta<n$,
we have
$$
\align
n  \int_M  F(q) &\log d(p,q)\,dV(q)\ \leq  \ \int_{\log F>-(n-\delta)\log d(p,q)}F(q) (-n\log d(p,q))\,dV(q)
\\&\qquad\qquad\qquad\qquad\qquad\qquad\qquad  +\ \int_{\log F\leq -(n-\delta)\log d(p,q)}F(q)
 (-n\log d(p,q))\,dV(q)\\
&\leq\ \frac{n}{n-\delta}\int_{M}F \log F\,dV\ +\ n\int_{M} (d(p,q))^{\delta-n} \log d(p,q)\,dV(q).
\endalign
$$
The second integral on the right converges, so
choosing $\delta$ small enough so $n/(n-\delta)< 1+\ep$, we get (2.2).
\medskip\medskip

\noindent{\bf Proof of Lemma 2.5}.  For each point $p_0\in M$ we there exist smooth
 coordinates $(y_1,\dots,y_n)$ mapping some open neighborhood of $p_0$ diffeomorphically onto an
open ball $B_{\Bbb R^n}(0,r)$.  Then
writing
$$
 g_{ij}=g\left( \partial_{y_i},\partial_{y_j}\right),
$$
we have
$$
dV\ =\  \sqrt{\left| \det g_{ij}\right|}\ dy_1\dots dy_n.
$$
Set
$$
x_i=y_i,\ \ i<n,\qquad\qquad\qquad x_n=\int_0^{y_n} \sqrt{\left|\det g(y_1,\dots,y_{n-1},t)\right|}\,dt.
$$
Then by the inverse function theorem, the map $y\to x$ is a diffeomorphism from some neighborhood
of zero to a neighborhood of zero, so $p\to x(y(p))$ defines smooth coordinates on
some neighborhood of $p_0$.  Moreover, by computing the Jacobian $|\partial y/\partial x|$ we see that
$$
dV\ =\ dx_1\dots dx_n.
$$
We now work in the coordinates $x$ and write
$$
g_{ij}=g(\partial_{x_i},\partial_{x_j}).
$$
By applying a linear transformation to $(x_1,\dots,x_n)$ if necessary we can assume that at $x=0$,
$g_{ij}=\delta_{ij}$.
Then there exists  $r>0$ such that for $|x|<r$
$$
\left(\sum_{i,j}| g_{ij}-\delta_{ij}|^2\right)^{1/2}\ <\ 1-e^{-\ep}.
$$
But then for $v=(v_1,\dots,v_n)$,
$$
e^{-\ep}|v|\ \leq\ g(v,v)^{1/2}\ \leq\ e^\ep|v|,
$$
and for any curve $\gamma$ in the coordinate ball $\{x:|x|<r\}$, we have that if
$|\gamma|$ is the Euclidean length of $\gamma$ and $L_g(\gamma)$ is
the length of $\gamma$ in the metric $g$, then
$$
e^{-\ep}|\gamma|\ \leq\ L_g(\gamma)\ \leq\ e^\ep |\gamma|.
$$
But then minimizing the middle term or the term on the right
over curves $\gamma$ joining $y$ to $z$ gives (2.5).
Now the ball $|x|<r$ contains some geodesic ball $B_g(p_0,\delta)$.

We can choose $\delta$ independent of $p_0$ by a simple compactness argument.
Indeed, by compactness there is a finite cover of balls $B_g(p_1,\delta_1),\,\dots,\, B_g(p_N,\delta_N)$
on which we have coordinates satisfying (2.4) and (2.5). But since this is an open cover,
we can find $\delta>0$
such that for each $p_0\in M$, $B_g(p_0,\delta)\subset B_g(p_j,\delta_j)$ for some $j\in \{1,2,\dots,N\}$.
\medskip\medskip

\noindent{\bf Proof of Proposition 2.7}. This is a special case of Proposition 2.12 proved below.
\medskip\medskip

\noindent{\bf Proof of Lemma 2.8} (a)$\Rightarrow$(b).  Just apply Cauchy-Schwarz to get.
$$
\multline
\frac1V\int_M Q A^{-1} R\,dV\ \leq\
\left(\frac1V\int_M Q A^{-1} Q\,dV\right)^{1/2}\left( \frac1V\int_M R A^{-1} R\,dV\right)^{1/2}
\\ \leq\ \frac12 \left(
\frac1V\int_M Q A^{-1} Q\,dV\ +\  \frac1V\int_M R A^{-1} R\,dV\right),
\endmultline
$$
and then apply (a).

\noindent (b)$\Rightarrow$(c).  Just apply (b) to $Q/\s_Q$ and $R/\s_R$.

\noindent (c)$\Rightarrow$(a).  (a) is just a special case of (c) when $Q=R=F$ and $V_F=V$.
\medskip\medskip

\noindent{\bf Proof of Proposition 2.9}.  Now
$\mu_V(\Bbb R^n)$ is bounded below (see [CL] or (2.11)).  Moreover, applying Lemma 2.4(a) with
$\ep=1$, we see
that there  exists $\delta>0$ such that if
and $p\in M$ and  $F\in (L\log L)^+(M)$ is supported in $B_g(p,\delta)$, then
$$
\mu(M,g,F)\ \geq\  \mu_{V_F}(\Bbb R^n)\ -\ V_F.
$$
Hence from the definition of $\mu(M,g,F)$ and the fact that $m$ is bounded, we get
a constant $C$ such that if $F\in (L\log L)^+(M)$ is supported in $B_g(p,\delta)$ with $V_F=V$, then
$$
\frac1{V}\int_M F A^{-1} F\,dV\ \leq\ \frac2{\gamma_n}\left( \int_M F \log F\,dV\right)
\ +\ C.
$$
This proves (2.12) when $F$ is supported in $B_g(p,\delta)$ for some $p\in M$, but we
want to remove this restriction on the support of $F$.
As in Lemma 2.8, we see that whenever  $Q,R\in (L\log L)^+(\Bbb R^n)$ with
$Q,R\in B_g(p,\delta)$ for some $p$, then
$$
\frac1V\int_M Q A^{-1} R\,dV\ \leq\ \frac1{\gamma_n}\left( \s_R\int_M Q\log \frac{Q}{\s_Q}\,dV
\ +\  \s_Q \int_M R\log \frac{R}{\s_R}\,dV\right)\ +\ C\,\s_Q\,\s_R.
$$
Choose closed sets $W_1,\dots,W_N$  which
cover $M$ such that the measure of $W_i\cap W_j $ equals zero if $i\neq j$. Suppose also
that the sets $W_j$ are sufficiently small that
 if $W_i\cap W_j\neq \emptyset$ then there exists $p$ with $W_i$ and $W_j$ contained in $B_g(p,\delta)$.
 We can choose $\ep>0$ such that if $W_i\cap W_j=\emptyset$ then the distance from $W_i$ to $W_j$ is
 at least $\ep$. Let $\chi_j$ denote
the characteristic function of $W_j$.  For $F\in (L\log L)^+(M)$,  set $F_j=\chi_j F$.
Set
$$
C_1\ =\ \sup_{d(p,q)>\ep} |G(p,q)|.
$$
Then
$$
\int_{M} F A^{-1} F\,dV
\ =\ \sum_{W_i\cap W_j=\emptyset}\int_M F_i A^{-1} F_j\,dV
\ +\ \sum_{W_i\cap W_j\neq\emptyset}\int_M F_i A^{-1} F_j\,dV.
$$
The first sum on the right is bounded by $C_1 V_F^2$.
The second sum on the right   is bounded by
$$
 \sum_{W_i\cap W_j\neq \emptyset} \biggl(
\frac{V_{F_i}}{\gamma_n}\int_{M} F_j\log \left( \frac{F_j}{\s_{F_j}}\right)\,dV
\ +\ \frac{V_{F_j}}{\gamma_n}\int_{M} F_i\log\left(\frac{ F_i}{\s_{F_i}}\right)\,dV\ +\
\frac{CV_{F_i}V_{F_j}}V\biggr). \tag 3.4
$$
However, since the function $t\to t\log t$ is convex, Jensen's inequality gives
$$
\int_M F_j \log\left( \frac{F_j}{\s_{F_j}}\right)\,dV\ \geq\ 0.
$$
Hence assuming $V_F=V$, (3.4) is bounded above by
$$
\align
&  \sum_{i,j}
\biggl(
\frac{V_{F_i}}{\gamma_n}\int_{M} F_j\log \left( \frac{F_j}{\s_{F_j}}\right)\,dV
\ +\ \frac{V_{F_j}}{\gamma_n}\int_{M} F_i\log\left(\frac{ F_i}{\s_{F_i}}\right)\,dV\ +\
\frac{CV_{F_i}V_{F_j}}V\biggr)\\
& =\ \frac{2V_{F}}{\gamma_n}\sum_j \int_{M} F_j\log \left( \frac{F_j}{\s_{F_j}}\right)\,dV\ +\
\frac{CV_{F}^2}V\\
&=\ \frac{2V}{\gamma_n}\int_{M} F\log F\,dV\ -\
\frac{2V}{\gamma_n}\sum_j V_{F_j}\log V_{F_j}\ +\ C'
\endalign
$$
where $C'$ depends only on $C$ and $V$.  Setting $C''=-\min_{t>0} t\log t$, we get
$$
\int_M F A^{-1} F\,dV\ \leq\
\frac{2V}{\gamma_n}\int_{M} F\log F\,dV\ +\ \frac{2NC''V}{\gamma_n} \ +\
C'\ +\ C_1 V^2.
$$
This completes the proof of Proposition 2.9.
\medskip\medskip

\noindent{\bf Proof of Lemma 2.11}.  Since the function $\Phi(t)$ grows faster at
infinity than the function $t$, we see that  functions $F_j$ are uniformly bounded in
$L^1(M)$. Hence by taking a subsequence we can
assume $F_j$ converges weakly to a measure $d\sigma$, that is for all $\phi\in C(M)$,
$$
\int_M F_j \phi\,dV\ \to \ \int_M \phi\,d\s\qquad \text{ as }j\to\infty. \tag 3.5
$$
But applying Tchebychev's inequality to the functions $F_j$, we see that for a
measurable set $U\subset M$,
$$
\int_U F_j\,dV\ \leq\  \int_{\{x\in U:F_j(x)\leq\la\}}F_j\,dV\ +\
\int_{\{x\in U:F_j(x)>\la\}}F_j\,dV\ \leq\ \la \int_U\,dV \ +\ \frac{S}{\Phi(\la)},
$$
and this shows that the limit $d\s$ is absolutely continuous with respect to the measure $dV$, and
hence equals $F\,dV$ for some function $F\in L^1(M)$, and (3.5) holds when $\phi$ is the
characteristic function of a measurable set.  Now for $0\leq m\leq N^2$,
$$
I^m_N\ =\ \left\{x: \frac{m}N\ \leq\ F\ \leq\ \frac{m+1}N\right\},\qquad\qquad
V^m_N=\int_{I^m_N}\,dV.
$$
Then by Jensen's inequality,
$$
\multline
\sum_{m=0}^{N^2} \Phi\left( \frac1{V^m_N }\int_{I^m_N } F \,dV\right)\ V^m_N
\ =\ \liminf_{j\to\infty} \sum_{m=0}^{N^2} \Phi\left(\frac1{V^m_N }\int_{I^m_N } F_j\,dV\right)
\ V^m_N \\ \leq\ \liminf_{j\to\infty} \sum_{m=0}^{N^2}\int_{I^m_N } \Phi(F_j)\,dV\ =\
\liminf_{j\to\infty}\int_{\{F\leq N^2\}} \Phi(F_j)\,dV
\endmultline
$$
However, this gives (2.18), because as $N\to\infty$ the left hand side of (3.8) converges to
the left hand side of (2.18) and the right hand side of (3.8) converges to the right hand side of
(2.18).
\medskip\medskip

\noindent{\bf Proof of Proposition 2.12}. From Lemma 2.3 we see that (2.24) holds
 if the supports of the functions $F_j$
are shrinking to a point.  In order
to prove equality for general concentrating sequences,
take $\delta_j\to 0$ with
$$
\int_{B_g(p_j,\delta_j)} F_j\,dV\ >\ V(1-\delta_j).
$$
Let $\chi$ be the characteristic function of $B(p_j,\delta_j)$ and set
$$
Q_j= \chi\,F_j,\qquad\qquad R_j =F_j-Q_j.
$$
Now we apply  Lemma 2.4. Since $\delta_j\to 0$ and $Q_j$ is supported in $B(p_j,\delta_j)$,
we can find a sequence $\ep_j>0$ with $\ep_j\to 0$ such that
$$
\mu(M,g,Q_j)\ \geq\ \mu_{V_{Q_j}}(\Bbb R^n)\ -\ \ep_j V_{Q_j}. \tag 3.6
$$
To simplify the notation, we fix $j$ and set $F=F_j$, $Q=Q_j$, $R=R_j$. Then
by Proposition 2.9 and Lemma 2.8, we have
$$
\int_M Q A^{-1} R\,dV\ \leq\ \frac1{\gamma_n}\left( V_R\int_M Q\log \frac{Q}{\s_Q}\,dV
\ +\  V_Q \int_M R\log \frac{R}{\s_R}\,dV\right)
\ +\ C\,V_Q\,V_R.\tag 3.7
$$
But since  $Q$ and $R$ have disjoint supports, one easily checks that
$$
\align
V&\ \mu(M,g,F)\ =\ V\int_M Fm\,dV\ +\ \frac{2 V_{F}}{\gamma_n} \int_M F\log F\,dV\ -\
\int_M F A^{-1}F\,dV\\  &=\ (V_Q+V_R)\int_M(Q+R)m\,dV\ +\ \frac{2(V_Q+V_R)}{\gamma_n}\int_M\left( Q\log Q
\ +\  R\log R\right)\,dV
\\&\qquad\qquad\qquad\qquad\qquad\qquad\qquad\qquad\qquad\qquad\qquad\qquad\qquad
 -\ \int_M (Q+R) A^{-1}(Q+R)\,dV\\
&=\ V_Q\ \mu(M,g,Q)\ +\ \frac{2V_R}{\gamma_n}\int_M R\log R\,dV\ -\ \int_M RA^{-1}R\,dV
\\ &\qquad\qquad +\ \frac{2V_Q}{\gamma_n}\int_M R\log R\,dV\ +\ \frac{2V_R}{\gamma_n}\int_M Q\log Q\,dV
\ -\ 2\int_M Q A^{-1} R\,dV\\
&\qquad\qquad +\ \int_M (V_Q R+ V_R F )m\,dV.
\endalign
$$
Then
applying  (3.7), we get
$$
\align
V\ \mu(M,g,F )\ &\geq\ V_Q\ \mu(M,g,Q)\ +\ \frac1{\gamma_n}\biggl(2V_R^2\log \s_R\ +\
2V_R V_Q(\log \s_R+\log \s_Q)\biggr)\\
&\qquad\qquad\qquad\qquad\qquad\qquad\qquad\qquad
-\ \frac{CV_R^2}V \ -\ \frac{2CV_R V_Q}V\ -\ V_R (V_Q +V)\max_M |m|\\
&\geq \ V_Q\ \mu(M,g,Q)\ +\ 2(V_R \log V_R )V_{F }\ -\ C'V_R ,
\endalign
$$
where $C'$ depends only on $(M,g)$. Allowing $j$ to vary and applying (3.6), we
have
$$
\mu(M,g,F_j)\ \geq\
\frac{V_{Q_j}}{V}\left(\mu_{V_{Q_j}}(\Bbb R^n)\ -\ \ep_j V_{Q_j}\right)
\ +\ 2(V_{R_j}\log V_{R_j} )\ -\ C'V_R.
$$
But the right hand side converges to $\mu_{V}(\Bbb R^n)$ as $j\to\infty$.
  This completes the proof of Proposition 2.12.
\medskip\medskip

\noindent{\bf Proof of Proposition 2.13}.  We  modify the proof of Proposition 2.9.
(Write $\delta'$ for the value of $\delta$ appearing in the proof of Proposition 2.9 to distinguish
it from the value $\delta$ in the statement of Proposition 2.13.)
When we choose the sets $U_1,\dots, U_N$ we place the
additional restriction that $U_j\subset B_g(p_j,\delta/2)$ for
some points $p_j\in M$.
Then choosing $V_j$ as above, we get that for fixed $j$,
$$
\sum_{\{i:V_i\cap V_j\neq \emptyset\} }V_{F_i}\ \leq\ \int_{B_g(p_j,\delta)} F\,dV\ \leq\ (1-\delta) V_F.
\tag 3.8
$$
Using this to bound (3.4), we get
$$
\multline
\int_M FA^{-1}F\,dV\ \leq\ \frac{ 2(1-\delta)V_F }{\gamma_n}\sum_j \int_M F_j \log \left(\frac{F_j}{\s_{F_j}}\right)\,dV
\ +\ \frac{CV}V\ +\ C_1V^2\\
\leq\ \frac{ 2(1-\delta)V}{\gamma_n}\int_M F\log F\,dV\ +\ C_2,
\endmultline
$$
where $C_2$ depends only on $(M,g)$ and $\delta$. This completes the proof of Proposition 2.13.

\medskip
\medskip

\noindent{\bf Proof of Lemma 2.14}.   We will assume that the reader is familiar with the standard theory
of elliptic pseudodifferential operators, in particular the construction of powers, see for
example [Se].    We work with  the spaces $\PDO^m(M)$ of classical
pseudodifferential operators of order $m\in \Bbb R$. This is a class of operators on
$C^\infty(M)$.   For local coordinates on $\Omega\subset M$,
an operator $B\in \PDO^m(M)$ acts on smooth functions $F$ supported in $\Omega$ by
$$
BF(x)\ =\ \frac1{(2\pi)^n}\int_{\Bbb R^n} \int_{\Omega} e^{i(x-y)\cdot \xi}\, b(x,\xi) F(y)\,dyd\xi,
$$
where $x,y$ are the coordinates of  points in $\Omega$ and
where the symbol $b(x,\xi)$ satisfies estimates
$$
\left| \partial_x^\al \partial_\xi^\beta b(x,\xi)\right|\ \leq\ C_{\al \beta} (1+|\xi|)^{m-|\beta|}.
$$
We say that a set of operators $B^{(\ep)}$ is {\it uniformly bounded in $\PDO^m (M)$} if the constants $C_{\al \beta}$
can be chosen independent of $\ep$. If $m<0$ then $B$ is an integral operator.  It's Schwartz
kernel $K(p,q)$ is a function on $M\times M$ which is given in local coordinates on $\Omega\times\Omega$
by
$$
K(x,y)\ =\ \frac1{(2\pi)^n}\int_{\Bbb R^n} e^{i(x-y)\cdot \xi}\, b(x,\xi)\,d\xi.
$$

For $r>0$, the operators $A^{-\ep}$ with $\ep\in [0,r]$ are
uniformly bounded in $\PDO^0(M)$.  The following standard Lemma is
not quite sufficient to prove Lemma 2.14.

\proclaim{Lemma 3.1}  If $B^{(\ep)}$ are uniformly bounded in $\PDO^m(M)$ and $m<k$, then the operators
$B^{(\ep)}$ are uniformly bounded from  $C^k(M)$ to $C(M)$.
\endproclaim

We start by proving Lemma 2.14 in the case $k=0$.
We want to apply Lemma 2.14 to $A^\ep$, but we will need to deal separately with the principal symbol.
To do this we will apply the following trivial lemma.

\proclaim{Lemma 3.2}  If for $F\in C(M)$ we set
$$
BF(p)\ =\ \int_M K(p,q) F(q)\,dV(q),
$$
where
$$
\sup_p \int_M |K(p,q)|\,dV(q)\ \leq\ C,
$$
then  $B$  is bounded on $C(M)$ and
$$
\|BF\|_\infty\ \leq\ C\|F\|_\infty.
$$
\endproclaim
Working in coordinates $(x_1,\dots,x_n)$,
we write   $g^{ij}$ for the components of $g^{-1}$. Then the principal symbol of $A^{-\ep}$
is $(g(x,\xi))^{-n\ep/2}$ where
$$
g(x,\xi)\ =\ \sum_{i,j}g^{ij}\xi_i\xi_j.
$$
Working in normal coordinates about the point $p_0$ then this symbol is just $|\xi|^{-n\ep}$.
The kernel of the operator corresponding to this symbol is given by taking the
inverse Fourier transform of $|\xi|^{-n\ep}$:
$$
 \frac1{(2\pi)^n} \int_{\Bbb R^n} e^{ix\cdot \xi} |\xi|^{-n\ep}\,d\xi
\ =\  C(\ep)\, |x|^{\ep-n}
\qquad\qquad
C(\ep)\ =\ (2\pi)^{-n/2} 2^{n(\frac12-\ep)}\frac{\Gamma((n(1-\ep)/2)}{\Gamma(\ep/2)}. \tag 3.9
$$
Notice that  $C(\ep)/\ep$ is bounded as $\ep\to 0$.
Here, $x$ are the normal coordinates on $M$ centered at  $p_0$, so $|x|$ measures the distance from $p_0$.
The upshot of this is that we can write
$$
A^{-\ep} \ =\ A^{-\ep}_0\ +\ B^{(\ep)},
$$
where
$$
A^{-\ep}_0 F(p)\ =\ C(\ep) \int_M   (d(p,q))^{\ep-n}\phi(d(p,q)) F(q)\,dV(q). \tag 3.10
$$
Here $\phi$ is a smooth cut off function, and the operators $B^{(\ep)}$ are uniformly bounded in the class
$\PDO^{-1}(M)$.
(There is a technical point here: we are writing $q$ in normal coordinates
around $p$ instead of using a fixed
coordinates chart for both $p$ and $q$.)

Applying Lemma 3.1 to $B^{(\ep)}$ we find that they are uniformly bounded on $C(M)$.
Applying Lemma 3.2 to (3.10) we find that $A^{-\ep}_0$ are uniformly bounded on $C(M)$.
Hence $A^{-\ep}$ are uniformly bounded on $C(M)$.

To show that $A^{-\ep}$ are uniformly bounded on $C^k(M)$ we just need to show that if
$D_k$ is a partial differential operator of order $k$ on $M$ then $D_k A^{-\ep}$ are uniformly bounded
from $C^k(M)$ to $C(M)$. But
$$
D_k A^{-\ep}\ =\ A^{-\ep} D_k \ +\ [D_k, A^{-\ep}] .
$$
Now $D_k$ is bounded from $C^k(M)$ to $C_M$, and $A^{-\ep}$ is uniformly
bounded from $C(M)$ to $C(M)$ so we have dealt with the term $A^{-\ep} D_k$.  The commutators
$[D_k,A^{-\ep}]$ are uniformly bounded in $\PDO^{k-1}(M)$, and so by Lemma 3.1 they are uniformly
bounded from $C^k(M)$ to $C(M)$.

\medskip\medskip
\noindent{\bf Appendix}.
\define\fp{{\operatorname{f.p.}}}

\medskip

Let $M$ be a closed,  manifold with metric $g$. Denote the volume
element by $dV$. Let $A$ be an operator of type $\Delta^{n/2}$ with
Robin mass $m$. In this appendix, we compute the difference between
the regularization
$$
\trace A^{-1}\ =\ \int_M m\,dV
$$
and the zeta function regularization defined as follows.  Let $\la_1\leq\la_2\leq\dots$
be the non-zero eigenvalues of $A$ and set
$$
Z(s)\ =\ \sum_j \la_j^{-s}.
$$
By Weyl's law,  $Z(s)$ converges for $\Re s>1$.  Now $Z(s)$
has an analytic continuation to a meromorphic function of $s\in \Bbb C$
with a simple pole at $s=1$. We define $\trace_\zeta A^{-1}$ to be the finite part of
$Z(s)$ at $s=1$, that is
$$
\trace_\zeta A^{-1}\ =\  Z(s)\bigl|^{\fp}_{s=1}\ :=\ \frac{d}{ds}\biggl|_{s=1}(s-1)Z(s).
$$
To compare these two regularizations of the trace of $A^{-1}$, let
$\phi_j$ be the eigenfunction of $A$ with eigenvalue $\la_j$ which is
normalized in $L^2(M)$.  The Schwartz kernel of $\Delta^{-s}$ is given by
$$
K(A^{-s},p,q)\ =\ \sum_j \la_j^{-s} \phi_j(p)\phi_j(q),
$$
where the sum converges in the distributional sense. Away from the diagonal $p= q$,
the function $K(A^{-s},p,q)$ is
entire in $s$ and smooth in $(s,p,q)$. At the diagonal, $K(A^{-s},p,p)$
has a analytic continuation to a meromorphic function and
$$
Z(s)\ =\ \int_M K(A^{-s} p,p)\,dV(p). \tag A.1
$$
Hence
$$
\trace A^{-1}\ -\ \trace_\zeta A^{-1}\ =\ \int_M \left( m(p)\ -\ K(A^{-s},p,p)\bigl|^\fp_{s=1}\right)\,dV(p)
\ =\ \int_M c(p)\,dV(p),
\tag A.2
$$
where $c(p)$ is the difference in two different ways of regularizing $K(A^{-1},p,p)$:
$$
c(p)\ =\
\lim_{q\to p} \left(K( A^{-1},p,q)\ +\ \frac{2n}{\gamma_n}\log  d(p,q)\right)\ -\ K(A^{-s},p,p)\bigl|^\fp_{s=1}.
\tag A.3
$$
The operator $A^{-s}$ is a pseudodifferential operator whose symbol expansion can be computed in
terms of the symbol of $A$.  In particular, working in normal coordinates at the point $p$, the principal
symbol of $A^{-s}$ is $|\xi|^{-ns}$.  If one is familiar with the construction of powers of
an elliptic operator, see for example [Se], it is not difficult to show that only the principal symbol
will contribute to the anomaly (A.3).  Indeed, if we define $K(s,x)$ to be the kernel of this symbol,
that is
$$
K(s,x)\ =\ \frac1{(2\pi)^n}\int_{|\xi|>1} e^{-ix\cdot \xi}|\xi|^{-ns}\,d\xi,
$$
and if $X(q)$ denotes the normal coordinate of $q$ centered at $p$, then
the function
$$
K(A^{-s},p,q)\ -\ K(s,X(q))
$$
is continuous in $(s,p,q)$ for $q$ in a neighborhood of $p$ and $s\geq 1$.
Hence from (A.3) we have
$$
c(p)\ =\ \lim_{x\to 0} \left(K(1,x)\ +\ \frac{2n}{\gamma_n} \log |x|\right)
\ -\ K(s,0)\bigl|^\fp_{s=1}.
\tag A.4
$$
Now  for $s<1$ the function $|\xi|^{-ns}$ defines a homogeneous distribution, and its inverse Fourier
transform can be computed by duality. As in (3.9),
$$
\frac1{(2\pi)^n}\int_{\Bbb R^n} e^{ix\cdot\xi}|\xi|^{-ns}\,d\xi\ =\
C(s) |x|^{n(s-1)},\qquad\qquad\qquad C(s)= (2\pi)^{-n/2}\,2^{n(\frac12-s)}
\frac{\Gamma(n(1-s)/2)}{\Gamma(ns/2)}.
$$
Hence computing in polar coordinates,
$$
\align
K(s,x)\ &=\ C(s)|x|^{n(s-1)}\ -\ \frac1{(2\pi)^n}\int_{|\xi|<1} e^{ix\cdot \xi} |\xi|^{-ns}\,d\xi
\\ &=\ C(s)|x|^{n(s-1)}\ -\  \frac{|S^{n-1}|}{(2\pi)^n}\int_{r<1}r^{-ns+n-1}\,dr\ +\ O(|x|^2)\\
&=\ C(s)|x|^{n(s-1)}\ -\ \frac{|S^{n-1}|}{n(2\pi)^n(1-s)}\ +\ O(|x|^2)\\
&=\ C(s)|x|^{n(s-1)}\ +\ \frac{2}{\gamma_n(s-1)}\ +\ O(|x|^2),
\endalign
$$
for small $|x|$. We see that
$$
K(s,0)\bigl|^\fp_{s=1}\ =\ 0.
$$
 We have
$$
c(p)\ =\ \lim_{x\to 0} \left( C(s)|x|^{n(s-1)}\ +\ \frac{2}{\gamma_n(s-1)}
     \ +\ \frac{2n}{\gamma_n} \log |x|\right)\ =\ C(s)\bigl|^\fp_{s=1}.
$$
From (A.2),
$$
\trace A^{-1}\ -\ \trace_\zeta A^{-1}\ =\  C(s)\bigl|^{\fp}_{s=1}\ V. \tag A.5
$$
Explicitly we compute
$$
C(s)\bigl|^\fp_{s=1}\ =\ \frac1{(4\pi)^{n/2} \Gamma(n/2)}\left( 2\log 2\ +\ \Gamma'(1)\ +\
\frac{\Gamma'(n/2)}{\Gamma(n/2)}\right). \tag A.6
$$

\medskip\medskip

\centerline{Acknowledgements}
\medskip

I would like to thank  S.Y. Alice Chang for information on GJMS
operators, and Niels Martin M{\o}ller for pointing out some
corrections in the exposition.
\medskip
\medskip

\centerline{References}
\medskip

\roster

\item"[Ad]" D. Adams:  A sharp inequality of J. Moser for higher order derivatives.
{\it Ann. of Math. (2)} {\bf 128} (1988),  385--398.

\item"[AH]" B. Ammann and E. Humbert: {Positive mass theorem for the Yamabe problem on spin manifolds}.
{\it Preprint}, (2003).

\item"[Au]" T. Aubin: Meilleures constantes dans le th\'{e}or\`{e}me  d'inclusion de Sobolev
et un th\'{e}or\`{e}me de Fredholm non lin\'{e}aire pour la transformation conforme de la courbure
scalaire. {\it J. Funct. Ana.} {\bf 32}, (1979) 148-174.

\item"[Be]" W. Beckner: Sharp Sobolev inequalities on the sphere and the Moser-Trudinger inequality.
{\it Annals of Math.} {\bf 138} (1993), 213-242.

\item"[Br]" T. Branson: An anomaly associated with 4-dimensional quantum gravity,
{\it Comm. Math. Physics}  {\bf 178} (1996), 301--309.

\item"[Br\O]" T. Branson and B. {\O}rsted,  Explicit functional determinants in four dimensions.
{\it Proc. Amer. Math. Soc.} {\bf 113} (1991),  669--682.

\item"[CL]" E. Carlen and M. Loss: Competing symmetries, the logarithmic HLS inequality
and Onofri's inequality on $S^n$.  {\it Geometric and Functional Analysis} {\bf 2} (1992)
90--104.

\item"[CC]" L. Carleson and S-Y. A.  Chang:
On the existence of an extremal function for an inequality of J. Moser.
{\it Bull. Sci. Math. (2)} {\bf 110} (1986), 113--127.

\item"[C]" Chang, Sun-Yung Alice: Conformal invariants and partial differential equations.
{\it Bull. Amer. Math. Soc.} {\bf  42} (2005),  365--393.

\item"[CQ]" S.-Y. A. Chang and J. Qing: Zeta functional determinants on manifolds with boundary.
{\it Math. Res. Lett.} {\bf 3} (1996), 1 -- 17.

\item"[CY]" S.-Y. A. Chang and P. Yang: Extremal metrics of zeta function determinants on 4-manifolds.
{\it Ann. of Math. (2)} {\bf 142} (1995), 171 -- 212.

\item"[DS1]"  P. Doyle and J. Steiner: Spectral invariants and playing hide and seek on surfaces.
{\it Preprint}.

\item"[DS2]"  P. Doyle and J. Steiner: Blowing bubbles on the torus.
{\it Preprint}.

\item"[FG]" C. Fefferman and  C. R. Graham, $Q$-curvature and Poincar\'e metrics.
{\it Math. Res. Lett.} {\bf 9} (2002),  139--151.

\item"[F]" L. Fontana, Sharp borderline Sobolev inequalities on compact Riemannian manifolds.
{\it Comment. Math. Helv.} {\bf 68} (1993),  415--454.

\item"[GJMS]" C. R. Graham, R.  Jenne,  L. Mason, and G. Sparling:
Conformally invariant powers of the Laplacian. I. Existence.
{\it J. London Math. Soc. (2)} {\bf 46} (1992),  557--565.

\item"[GZ]" C. R. Graham and M. Zworski: Scattering matrix in conformal geometry.
{\it Invent. Math.} {\bf 152} (2003),  89--118.

\item"[Gu]" M. J. Gursky, The principal eigenvalue of a conformally
invariant differential operator, with an application to semilinear
elliptic PDE. {\it Comm. Math. Phys. } {\bf 207}  (1999),  no. 1,
131--143.

\item"[H]"  L. Habermann:  Riemannian metrics of constant mass and
moduli spaces of conformal structures.
{\it  Lecture Notes in Mathematics}, {\bf 1743}. Springer-Verlag, Berlin, 2000.

\item"[HZ]" A. Hassell and S. Zelditch:  Determinants of Laplacians in exterior domains.
{\it Internat. Math. Res. Notices}  {\bf 18} (1999), 971--1004.

\item"[L]" E. Lieb: Sharp constants in the Hardy-Littlewood-Sobolev and related inequalities.
{\it Annals of Math.}, {\bf 118} (1983), 349--374.

\item"[Mor1]"  C. Morpurgo:  The logarithmic Hardy-Littlewood-Sobolev inequality
and extremals of zeta functions on $S^n$.
{\it Geom. Funct.} Anal. {\bf 6} (1996),  146--171.

\item"[Mor2]" C. Morpurgo:
 Sharp inequalities for functional integrals and traces of conformally invariant operators.
 {\it Duke Math. J.} {\bf 114} (2002),  477--553.

\item"[Mos]" J. Moser: A sharp form of an inequality by N. Trudinger.
{\it Indiana Math. J.} {\bf 20} (1971),
1077--1092.

\item"[Ok*]" G. Okikiolu: {\it Aspects of the Theory of Bounded Integral Operators in $L^p$ Spaces}.
Academic Press. N.Y., 1971.

\item"[Ok1]" K. Okikiolu: Critical metrics for the determinant of the Laplacian in odd dimensions.
{\it  Ann. of Math. (2)}  {\bf 153} (2001), 471 -- 531.

\item"[Ok2]" K. Okikiolu:   Hessians of spectral zeta functions. {\it Duke Math. J.}
{\bf  124} (2004), 517--570.

\item"[Ok3]" K. Okikiolu: A negative mass theorem for  the
$2$-Torus. {\it Preprint}.

\item"[OW]" K. Okikiolu and C. Wang:  Hessian of the zeta function for the Laplacian on forms.
{\it Forum Math.} {\bf  17} (2005),  105--131.

\item"[On]" Onofri: E. On the positivity of the effective action in a theory of random surfaces.
{\it Comm. Math. Phys.} {\bf  86} (1982),  321--326.

\item"[OPS1]"  B. Osgood, R. Phillips and P. Sarnak: Extremals of determinants of Laplacians.
{\it J. Funct. Anal.} {\bf 80} (1988),  148--211.

\item"[OPS2]"  B. Osgood, R. Phillips and  P. Sarnak:
 Moduli space, heights and isospectral sets of plane domains.
{\it Ann. of Math. (2)} {\bf  129} (1989),  293--362.

\item"[R]" K. Richardson: Critical points of the determinant of the Laplace operator.
{\it J. Funct. Anal.} {\bf 122} (1994), 52 -- 83.

\item"[Sc]" R. Schoen: Conformal deformation of a Riemannian metric to constant scalar curvature.
{\it J. Diff. Geom.}, {\bf 20} (1984), 479-495.

\item"[SY]" R. Schoen and S.-T. Yau: On the proof of the positive mass conjecture in general relativity.
{\it Comm. Math. Phys.}, {\bf 65} (1979) 45--76.

\item"[Se]" R. Seeley: Complex powers of an elliptic operator. 1967 Singular Integrals
{\it Proc. Sympos. Pure Math., Chicago, Ill., 1966} pp. 288--307 Amer. Math. Soc., Providence, R.I.

\item"[St1]"  J. Steiner: {\it Green's Functions, Spectral Invariants, and a Positive Mass on Spheres}.
Ph. D. Dissertation, University of California San Diego, June 2003.

\item"[St2]" J. Steiner:  A geometrical mass and its extremal properties for metrics on $S\sp 2$.
{\it Duke Math. J.}  {\bf 129} (2005), 63--86.

\item"[T]"  M. Taylor: {\it Partial differential equations. II. Qualitative studies of linear equations.}
Applied Mathematical Sciences, {\bf 116}. Springer-Verlag, New York, 1996.

\item"[Y]" H. Yamabe: On a deformation of Riemannian structures on compact manifolds.
{\it Osaka Math. J.}, {\bf 12} (1960), 21-37.

\endroster
\medskip\medskip

\centerline{Kate Okikiolu}

\centerline{University of California, San Diego}

\centerline{okikiolu\@math.ucsd.edu}

\end